\documentclass[12pt,oneside,mathscr]{amsart}

\usepackage{amsthm,amsmath,amssymb,euscript,amscd,a4wide}

\usepackage[all]{xy}
\usepackage{multirow,bigstrut}

\usepackage{tikz,mathrsfs}
\usetikzlibrary{arrows,decorations.pathmorphing,decorations.pathreplacing,positioning,shapes.geometric,shapes.misc,decorations.markings,decorations.fractals,calc,patterns}
\usetikzlibrary{tikzmark}
\usetikzlibrary{decorations.text,calc,arrows.meta}

\newtheorem{thm}{Theorem}[section]
\newtheorem{cor}[thm]{Corollary}
\newtheorem{lemma}[thm]{Lemma}
\newtheorem{prop}[thm]{Proposition}

\theoremstyle{definition}
\newtheorem{dfn}[thm]{Definition}

\newtheorem{example}[thm]{Example}

\newtheorem{remark}[thm]{Remark}

\numberwithin{equation}{section}

\let\strong=\textbf
\newcommand{\C}{\mathcal{C}}

\newcommand{\RR}{\mathbb{R}}

\newcommand{\A}{\mathcal{A}}

\newcommand{\compo}{\stackrel{}{\raise 2pt\hbox{$\scriptscriptstyle\circ$}}}
\newcommand{\sra}[1]{\stackrel{\scriptstyle #1}{\rightarrow}}
\newcommand{\lra}[1]{\stackrel{\scriptstyle #1}{\longrightarrow}}
\renewcommand{\leq}{\leqslant}
\renewcommand{\geq}{\geqslant}
\newcommand{\arr}{\operatorname{arr}}
\newcommand{\ob}{\operatorname{ob}}
\newcommand{\dom}{\operatorname{dom}}
\newcommand{\codom}{\operatorname{codom}}

\newcommand{\dsimplex}{\overrightarrow{\Delta}}
\newcommand{\drsimplex}{\overrightarrow{R}}
\newcommand{\polygon}{\overrightarrow{P}}

\newcommand{\smat}[1]{\left(\begin{smallmatrix}#1\end{smallmatrix}\right)}
\newcommand{\mat}[1]{\begin{pmatrix}#1\end{pmatrix}}

\begin{document}

\title{Fully Triangulated Categories}
\author[A Maciocia]{Antony Maciocia}
\author[H Alzubaidi]{Husniyah Alzubaidi}
\address{Department of Mathematics and Statistics\\
The University of Edinburgh\\
The King's Buildings\\ Mayfield Road\\ Edinburgh, EH9 3JZ.\\}
\email{A.Maciocia@.ed.ac.uk}
\thanks{}

\date{\today}
\begin{abstract}
We modify the axioms of triangulated categories to include both higher triangles and distinguished maps of 
higher triangles. The distinguished maps are specializations of Neeman's ``good'' maps of $2$-triangles. The axioms
both simplify Neeman's axioms in his 1991 paper and generalize them to higher triangles in the way proposed by Balmer et al. We provide a geometric formulation via directed truncated simplices to enable a more concrete approach. The axioms are modelled by homotopy and derived categories. We look at a number of key theorems including the fact that sums of maps of $2$-triangles
are distinguished if and only if the maps are distinguished and a strong version of the 3x3 lemma. These illustrate some key
proof methods. We also show that maps of faces of distinguished triangles are distinguished. We construct some useful distinguished 5-triangles.
\end{abstract}

\maketitle

\section{Introduction}{}
The notion of triangulated category has played a central role across many areas of mathematics. While they
have a significant drawback in that the cone of a morphism is not natural, they still provide a useful abstract framework
to describe a range of mathematical phenomena. 

The established expression of triangulated categories is via axioms which were introduced by Verdier in his thesis and reproduced in \cite{VerdierAst}. These
axioms, often abbreviated to TR1-TR4, have been used more or less unchanged since then. An alternative axiom system was proposed by Neeman and expounded in his textbook \cite{NeemanBook}. This involved replacing the octahedral axiom TR4 with a stronger version of TR3. The usual axiom TR3 says that a partial map $(f,g,-)$ if two distinguished triangles can be completed to a ``good'' map of the triangles. In the book he shows how to deduce TR4. Harder, is the converse, that the stronger TR3 follows from TR4. This he did in his 1991 paper \cite{Neeman1991}. In that paper he studied the key properties of good maps and proposed another axiom system more suited to the homotopy and derived categories. But these are
rather complicated based on equivalence classes of objects. Neeman's notion of good map is one whose mapping cone is also distinguished. We propose to replace this with a more abstract notion of distinguished map. We will show that distinguished maps are always good but the converse need not hold. Neeman observed that modifying a good map via what he calls a ``lightning strike'' (we also adopt this terminology) preserves the fact that it is good. But he speculates that the difference of two good maps need not be a lightning strike (even up to isomorphism). We show that distinguished maps form a single lightning strike equivalence class. 

Neeman's axioms also show that there is a link between suitable maps of triangles and higher dimensional diagrams (in his case, octahedra). We generalise this idea by using higher triangles in all dimensions. We will call these $n$-triangles, where $n$ is the dimension. These have appeared as diagrams since the early days of triangulated categories, for example, in
\cite[Remark 1.1.4]{BBD}. They have also existed before triangulated categories in Homotopy Theory in the form of Postnikov systems, see \cite[IV, Ex 1]{GM}, where the authors also note that these are special diagrams in a triangulated category which they call ``hypersimplices''. There are other notions of higher dimensional diagrams in additive categories such as $n$-angles (see \cite{BT1}) and distinguished $n$-triangles of Balmer \cite{Balmer}, Kunzer \cite{K07} and Maltsiniotis \cite{M}.  A sketch that the derived category satisfies the axioms was given by Maltsiniotis in a preprint \cite{M} and completed by Groth and \v{S}\v{t}ov\'{\i}\v{c}e \cite{GS}. But these definitions do not seem to work for the homotopy category without the further refinement of distinguished maps. One of the aims of this work, and indeed, one of the main motivations for the work, is to propose a new set of axioms and we call additive categories with a compatible shift functor \strong{fully triangulated} if they are satisfied. With this refinement it is possible to show that the homotopy category of an additive category is fully triangulated in a direct way in line with the traditional proofs that the homotopy category is triangulated. We describe the proof of this but leave the technical details for \cite{Hus}. The thesis also verifies that the axioms allow for a direct construction of the functor $D^b(\mathcal{A})\to \mathcal{T}$ extending the inclusion of an abelian category $\mathcal{A}$ in a fully triangulated category $\mathcal{T}$ as Neeman does in \cite{Neeman1991}.

The notion of distinguished $n$-triangles is a key feature of the axioms used by Balmer, Maltsiniotis et al. While our distinguished $n$-triangles are essentially the same, there is a key link to the new notion of distinguished maps. To motivate this consider an octahedral diagram arising from the usual octahedral axiom. If we pick two composable maps in this octahedron then the octahedral axiom allows us to complete this to octahedra. The non-uniqueness of cones on the three maps means that there is a potentially large class of different octahedra based on these two maps (and their composite). It is natural to ask if they are isomorphic as would be true for $2$-triangles determined by the the same single morphism. We can construct isomorphisms given by the four $2$-triangles in the octahedron but there is no reason why these can be assembled into a map of the two octahedra. In fact, even if we ask for the commuting squares of the octahedra to be homotopy squares (in other words the maps of the 2-triangle faces are good in the sense of Neeman) it still does not follow that there is an isomorphism. Indeed, an explicit counterexample has been constructed in \cite{K09}. This is where distinguished $3$-triangles are needed. But then the maps of the faces of these special $3$-triangles are not general maps. In fact, they are distinguished in our new sense. It then follows, for example, that a distinguished $3$-triangle has homotopy squares although the converse need not hold.  

In other treatments of higher triangulations it is usual to state the axioms indexed by the dimension $n$. We do the same but do not use the notion of $n$-triangulated categories to mean that the axioms hold up to and including $n$. Instead we require the axioms to hold for all $n$. Our axioms for $n=2$ imply that the category is triangulated in the usual sense essentially in the same way that Neeman's axioms do. However, to show properties of distinguished maps even in dimension $2$, we need higher triangles, and notably $5$-triangles. So it makes sense to require these higher triangle axioms to hold even when we are only studying small values of $n$. One key result is to show that the direct sum of two maps of $2$-triangles is distinguished if and only if each of the factors are distinguished.
Finally, we use the 3x3 lemma as a test for the strength of the axioms, just as Neeman does for his notion of good maps, and we show the strongest form of the lemma namely that we can choose the maps of the triangles extended from a commuting square and then these can always be completed to a 3x3 diagram so long as the maps of triangles are distinguished (see Theorem \ref{t:3x3} in the final section). A key result we use to establish this is that summands of distinguished maps of $2$-triangles are distinguished.

Higher triangles are not just a theoretical curiosity: we would argue that they are a useful computation tool. We have illustrated this in this paper but further examples can be found in \cite{MJAM} where they are used in the context of the derived category of coherent sheaves to help establish technical results about Bridgeland stability.

We start by providing an accessible and fairly elementary geometric account of the data needed for higher triangles (or $n$-triangles).
\section{Higher triangles}
There are several attempts in the literature to define higher dimensional triangle diagrams in
categories with shift. In this section we recall these definitions and set up the notation
needed for the rest of the paper. We start with diagrams in Euclidean space.
\subsection{Simplices and Rectified Simplices}
\begin{dfn}
For any $n\geq0$, let $\Delta_n$ be the standard $n$-simplex in $\RR^{n+1}$ whose
vertices are $a_i=(0,\ldots,0,1,0,\ldots,0)$ for $0\leq i\leq n$ and
edges $e_{ij}$, $0\leq i<j\leq n$. We let $\overrightarrow{\Delta}_n$
be the directed graph skeleton of $\Delta_n$ with edges $e_{i,j}$ (so the direction is
$a_i\to a_j$ for $i<j$). So we can also view $\dsimplex_n$ as an $n+1$-gon $a_0\to a_1\to a_2\to\cdots a_n$ and
$a_0\to a_n$ together with all ``composites'' inside the polygon. For
any integer $m\in\mathbb{Z}$ we can rotate $\dsimplex_n$ to form a
different directed graph on the 1-skeleton of $\Delta_n$. We replace
$a_0\to a_n$ with the opposite and $a_{m-1}\to a_m$ with its
opposite on the polygon and then complete the graph with the obvious
composites. We denote the resulting graph by $\dsimplex_n^m$ and call
it the \emph{$m$-rotated $n$-simplex}. Note
that $\dsimplex_n^m=\dsimplex_n^{m+k(n+1)}$ for any integer $k$ but we
will see that it is more natural to lift the action to the
integers. We shall denote the directed polygons by $\polygon_n^m$ and underline the
observation that $\polygon_n^m$ determines $\dsimplex_n^m$.
\end{dfn}
\begin{example} $\dsimplex_0$ is a point, $\dsimplex_1$ is a single directed edge, $\dsimplex_2$ is
  an oriented completely asymmetric triangle and
  $\dsimplex_3$ is an oriented tetrahedron with a source vertex and a sink vertex.
 \end{example}
There are well defined face maps $f_i:\Delta_{n-1}\to \Delta_n$ by
inserting a $0$ in the $i$th coordinate. There are also degeneracy
maps $d_i:\Delta_n\to \Delta_{n-1}$ given by adding together the $i$th
and $(i+1)$st coordinates. Note that the maps restrict to the directed
simplex and define maps of directed graphs.

In the case $n=2$, there are exactly two orientations of $\Delta_2$. One is $\dsimplex_2$ and
the other is the rotationally symmetric one. We would like to view this not as a simplex but as
a rectified simplex.
\[\begin{tikzpicture}[yscale=0.65,xscale=0.5]
        \begin{scope}[thick,decoration={
                markings,
                mark=at position 0.5 with {\arrow{>}}},
            ] 
        \draw[postaction={decorate}] (0,0)--(4,0); 
           \draw[postaction={decorate}] (4,0)--(8,0); 
        \draw[postaction={decorate}] (4,0)--(2,2);
        \draw[postaction={decorate}] (0,0)--(2,2);
        \draw [postaction={decorate}](2,2)--(4,4);
        \draw [postaction={decorate}](6,2)--(4,0);
        \draw [postaction={decorate}](2,2)--(6,2);
        \draw[postaction={decorate}] (6,2)--(4,4);
        \draw[postaction={decorate}] (8,0)--(6,2);
    \end{scope}
        \node[left] at (0,0) {$a_{0}$};
        \node[right] at (8,0) {$a_1$};
        \node[above] at (4,4) {$a_2$};
        \node[below] at (4,0) {$a_{0,1}$};
        \node[left] at (2,2) {$a_{0,2}$};
        \node[right] at (6,2) {$a_{1,2}$};
        \node[above] at (4,2) {$\scriptstyle e_{2,0,1}$};
        \node[left] at (3,0.8) {$\scriptstyle e_{0,1,2}$};
        \node[right] at (5,0.7) {$\scriptstyle e_{1,0,2}$};
       \node[left] at (5,0.9) {$\scriptstyle [1]$};
    \end{tikzpicture}\]
\begin{dfn}
A \emph{rectified $n$-simplex} $R_n$ is formed from $2\Delta_n$ by
forming a new vertex midway along each edge of $\Delta_n$ and cutting each
original vertex off by a hyperplane through these midpoints. So the vertices of
$R_n$ are at points $a_{i,j}=(0,\ldots,1,\ldots,1,\ldots,0)$ with
$0\leq i< j\leq 0$ (corresponding to each $e_{i,j}$ of $\Delta_n$).
The edges can be listed as $e_{i,j,k}$ for $0\leq i\leq n$  and $0\leq
j<k\leq n$ with $j\neq i\neq k$. So $e_{i,j,k}$ presents the edge
from $a_{i,j}$ or $a_{j,i}$ to $a_{i,k}$ or $a_{k,i}$. The corresponding directed graph is
the 1-skeleton $\drsimplex_n$ directed by $e_{i,j,k}$ if $i<j$ or $k<i$
and  the opposite of $e_{i,j,k}$ if $j<i<k$. We will want to keep track
of which arrows are which so we will write the first type as
$a_{i,j}\to a_{i,k}$ or $a_{j,i}\to a_{k,i}$ and the latter as
$a_{i,k}\rightsquigarrow a_{j,i}$ or $a_{i,k}\lra{[1]} a_{j,i}$
\end{dfn}
\begin{example}
$\drsimplex_0$ is empty, $\drsimplex_1$ is a point, $\drsimplex_2$ is the cyclic directed
  triangle as illustrated in the diagram above, and $\drsimplex_3$ is a directed octahedron as indicated in the diagram.
$$\xymatrix@R=1.8pc@C=2pc{&a_{1,3}\ar[dr]\ar^(0.6){[1]}[ddl]|!{[dr];[dl]}\hole |!{[ddd];[dl]}\hole \\
a_{1,2}\ar[ur]\ar_{[1]}[d]&&a_{2,3}\ar^{[1]}[ll]\ar_(0.4){[1]}[ddl]\\
a_{0,1}\ar[rr]| !{[u];[dr]}\hole |!{[urr];[dr]}\hole \ar[dr]&&a_{0,3}\ar[u]\ar[uul]|!{[u];[dl]}\hole |!{[u];[ull]}\hole\\
&a_{0,2}\ar[uul] \ar[ur]}$$
\end{example}

The following is easy to prove (induct on $i$ for (3)).
\begin{prop}
A truncated $k$-simplex has 
\begin{enumerate}
\item $\displaystyle \binom{n+1}{2}$ vertices,
\item $\displaystyle \binom{n+1}{2}(n-1)$ edges,
\item $\displaystyle \binom{n+1}{i+1}(n+1-i)$ $i$-faces.
\item $n+1$ $(n-1)$-simplex faces and $n+1$ truncated $(n-1)$-simplex faces.
\end{enumerate} 
\end{prop}
This time we have two types of face maps $f^\sigma_i:\Delta_{n-1}\to R_n$
and $f^\pi_i:R_{n-1}\to R_n$. The first is given by inserting a $1$ in
position $i$ and the second is given by inserting $0$ in position
$i$. $f^\sigma_0$ and $f^\pi_0$ extend to give maps of directed graphs $\dsimplex_{n-1}\to\drsimplex_{n}$ and
$f_0^\sigma:\drsimplex_{n-1}^0\to\drsimplex_n$.

Observe that the vertices of $f^\sigma_0(\dsimplex_{n-1})$
are $a_{0i}$ for $1\leq i\leq n$. Whereas the vertices of
$f^r_0(\drsimplex_{n-1})$ are $a_{ij}$ for $0<i<j\leq n$. So we can represent these in
2-dimensions as in figure 2.

\begin{figure}[h]
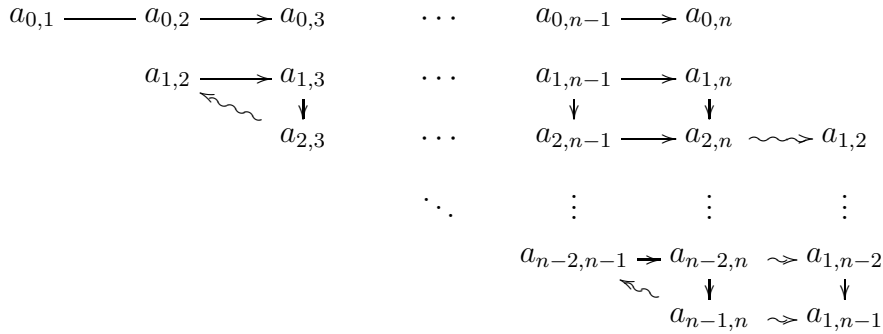
\label{fig:face}\xygraph{ 
!~-{@{-}@[|(2]}  !{0;/r1.8cm/:,p+/u0.8cm/::}
*+{a_{0,1}}="a1"
& *+{a_{0,2}}="a2"
& *+{a_{0,3}}="a3" &*+{\cdots} & *+{a_{0,n-1}}="an1" & *+{a_{0,n}}="an" \\
& *+{a_{1,2}}="a12" & *+{a_{1,3}}="a13" &*+{\cdots} & *+{a_{1,n-1}}="a1n1" &
*+{a_{1,n}}="a1n"   \\
&& *+{a_{2,3}}="a23" & \cdots & *+{a_{2,n-1}}="a2n1" & *+{a_{2,n}}="a2n" &
 *+{a_{1,2}}="a121" \\
&&&\ddots&\vdots&\vdots&\vdots\\
&&&& *+{a_{n-2,n-1}}="an2n1"
&*+{a_{n-2,n}}="an2n"& *+{a_{1,n-2}}="a1n21"\\
&&&&&*+{a_{n-1,n}}="an1n" &  *+{a_{1,n-1}}="a1n11"\\
"a1":@{-}"a2" "a2":"a3" "an1":"an"
"a12":"a13" "a1n1":"a1n" "a13":"a23" "a1n1":"a2n1" "a1n":"a2n" "a2n1":"a2n"
"an2n1":"an2n" "an2n":"an1n" 
"a23":@{~>}"a12" "an1n":@{~>}"an2n1"
"a2n":@{~>} "a121" "an2n":@{~>} "a1n21"
"an1n":@{~>} "a1n11" "a1n21":"a1n11"
  } \caption{The $0^{\text{th}}$ simplex and $0^{\text{th}}$ rectified simplex}
  \end{figure}
The other edges are given as composites vertically or horizontally. The
remaining edges (other than those given as composites of the ones
given) then make up the $\drsimplex_n$.
\begin{figure}[h]\label{fig:balmer}\xygraph{ 
!~-{@{-}@[|(2)]} !{0;/r1.8cm/:,p+/u0.8cm/::}
*+{a_{0,1}}="a1"
& *+{a_{0,2}}="a2"
& *+{a_{0,3}}="a3" &*+{\cdots} & *+{a_{0,n-1}}="an1" & *+{a_{0,n}}="an" \\
& *+{a_{1,2}}="a12" & *+{a_{1,3}}="a13" &*+{\cdots} & *+{a_{1,n-1}}="a1n1" &
*+{a_{1,n}}="a1n" & *+{a_{0,1}}="a011" \\
&& *+{a_{2,3}}="a23" & \cdots & *+{a_{2,n-1}}="a2n1" & *+{a_{2,n}}="a2n" &
*+{a_{0,2}}="a021" \\
&&&\ddots&\vdots&\vdots&\vdots\\
&&&& *+{a_{n-2,n-1}}="an2n1" &*+{a_{n-2,n}}="an2n"&*+{a_{0,n-2}}="a0n21"\\
&&&&&*+{a_{n-1,n}}="an1n" & *+{a_{0,n-1}}="a0n11"\\
&&&&&&*{a_{0,n}}="a0n1"
"a1":"a2" "a2":"a3" "an1":"an"
"a12":"a13" "a1n1":"a1n" "a13":"a23" "a1n1":"a2n1" "a1n":"a2n" "a2n1":"a2n"
"a011":"a021"
"an2n":"an1n" "a0n21":"a0n11" "a0n11":"a0n1"
"a23":@{~>}"a12" "an1n":@{~>}"an2n1"
"a2":"a12" "a3":"a13" "an1":"a1n1" "an":"a1n" "an2n1":"an2n"
"a1n":@{~>}"a011" "a2n":@{~>}"a021" "an2n":@{~>}"a0n21"
"an1n":@{~>}"a0n11" "a12":@{~>}"a1"
  } \caption{ Balmer diagram of $\drsimplex_n$}
\end{figure}

These diagrams were introduced by Paul Balmer \cite{Balmer} and have a number of useful features
which we will use extensively in what follows. One is that each $\dsimplex_{n-1}$ face is given
by the horizontal and vertical lines (extended to the right and above). The $\drsimplex_n$ faces
are a little harder to spot but the $\drsimplex_{n-1}$ faces are obtained by deleting the
row with row label $i$ and the column with the column number $i$. The $F_n$ face and $F_0$ faces
are easy to spot as indicated in figure 1.

\begin{prop}
  In an $\drsimplex_n$
  \begin{enumerate}
\item Each vertex appears exactly once in $f^\sigma_0(\dsimplex_{n-1})\cup
f^r_0(\drsimplex_{n-1})$.
\item Each edge of  arises exactly once in $\bigcup_i f^\sigma_i(\dsimplex_{n-1}^i)$.
\end{enumerate}
\end{prop} 

Note also that the (directed) valence of each vertex of $\drsimplex_n$ is always zero.
We also see from these diagrams that there are face maps
$f^\sigma_i:\dsimplex_{n-1}^i\to\drsimplex_n$. We see also that there
are rotated versions of $\drsimplex_n$ and face maps
$f^\pi_i:\drsimplex^i_{n-1}\to\drsimplex_n$. There is also a
canonically defined cycle:
\[a_{0,1}\to a_{0,n}\to a_{n-1,n}\rightsquigarrow
a_{n-2,n-1}\rightsquigarrow\cdots\rightsquigarrow a_{2,3}\rightsquigarrow
a_{1,2}\rightsquigarrow a_{0,1}\]
There is also a corresponding cycle of maps of faces $F_n\to F_{n-1}\to \cdots\to F_0\to
\tau^{n-1} F_n$, where $\tau$ is the rotation map defined in detail in the next section.

\subsection{Simplicial and rectified simplicial diagrams in categories with shift functors}
For any
category with shift endofunctor $(\C,[1])$, an \strong{$R_n$-diagram} or \strong{rectified simplex diagram}
is a diagram (graph) which is isomorphic as a graph to $\drsimplex_n$ whose edges are arrows or arrows
shifted by $[1]$ for the edge $e_{ijk}$ and $j<i<k$. We require that
each of the simplex faces are commuting diagrams and we also require
that the squares in the diagram commute as well.

We can also consider $\Delta^m_k$-diagrams which are simply commuting
diagrams in which the single reversed edge of the polygon has degree
$1$ when $m\ncong 0\bmod k+1$. Maps between them are just commuting
diagrams of maps between corresponding objects. We denote the category of these by
$\Sigma_k^m$ and abbreviate $\Sigma_k^0$ by $\Sigma_k$. Note that the
edges of this consists of actual maps.

The set of $R_k$-diagrams also form categories in the
obvious way (so that a map between diagrams is a collection of maps
between each objects for a given vertex which commutes with the edge
maps). We denote these categories by $\Pi_k$. Then the face maps
induce functors $F_m^\sigma:\Pi_k\to\Sigma_{k-1}^m$ and $F_m^\pi:\Pi_k\to\Pi_{k-1}$.
Similarly there are degeneracy maps $\Pi_{k-1}\to\Pi_k$. We can think of these as extending an $n-1$-triangle by zero maps or equality maps.
There is a good example of this in the proof of Proposition \ref{p:Nddash} below illustrating iterating the degeneracy map from a $2$-triangle.

There are also the \strong{rotation} $\tau$ and the \strong{full rotation} $\sigma$. Given an
$n$-triangle $A$ with base $F_0^\sigma(A)$, then $\tau(A)$ has $0$th base given by
$F_n^\sigma(A)$ and all maps are the same except for $e_{i,j,k}$ with $j<i<k$ in which case the
maps is $-e_{i,j,k}$. The full rotation $\sigma(A)$ has objects $a_{i,j}[1]$ in position $(i,j)$
and maps $-e_{i,j,k}$ for $j<i<k$. Then $\tau^{n+1}=\sigma^{n-1}$ and $\sigma^2=[2]$. 
Note that for $n=2$, $\sigma=\tau^3$.

\section{Fully triangulated categories}
We can now define the notion of fully triangulated category. Essentially it is a choice of
special $R_k$-diagrams and certain maps between them which we call distinguished.
\begin{dfn}
  We say that an category with endofunctor $(\C,[1])$ is
\strong{fully triangulated} if for each integer $n\geq2$ there is a
subcategory $D_n\subset \Pi_n$ consisting of 
\strong{distinguished $n$-triangles} and \strong{distinguished morphisms} of $n$-triangles
satisfying the following axioms for all $n\geq2$:
\begin{enumerate} 
\item\label{a:isom} [\strong{Isomorphism}] For any $A,B\in \ob(D_n)$, any isomorphism $A\to B$ is distinguished.
\item\label{a:rot} [\strong{Rotation}] For any $A\in \ob(D_n)$, the rotation $\tau(A)$ and the full rotation $\sigma(A)$ are in $\ob(D_n)$.  
\item\label{a:faces} [\strong{Faces}] For $n>2$, the image of a distinguished map of $n$-triangles under a face map
  and the image of a distinguished map of $n-1$-triangles under a degeneracy map is
  a distinguished map. 
\item\label{a:bases} [\strong{Bases}] Every $n-1$-simplex diagram $\Sigma$ extends to a distinguished $n$-triangle and every map of $n-1$-simplex diagrams $f:\Sigma\to\Sigma'$ extends to a distinguished
  map of any distinguished $n$-triangles extending the bases. 
\item\label{a:ls} [\strong{Lightning Strike}] If $(u,v,w,f,g,h,u',v',w')$ and $(u,v,w,f,g,h',u',v',w')$
  are two distinguished maps of the same $2$-triangles then there is a map
  $\tau:\codom(w)\to\dom(v')$ such that $h'-h=v'\tau w$.   
\end{enumerate}
\end{dfn}

In what follows $\C$ will always be a fully triangulated category and we leave the rest of the structure understood.
\begin{example}
Let $\A$ be any additive category. Then the homotopy categories $K^*(\A)$ (where the ${}^*$ is either empty, $+$, $-$ or $b$) are fully triangulated by defining
distinguished $n$-triangles to be $n$-triangular diagrams isomorphic to certain standard ones constructed via the
mapping cones of maps. Distinguished morphisms of $2$-triangles are maps of $2$-triangles isomorphic to certain maps of the
standard $2$-triangles: 
\begin{equation*}
\xymatrix@=1.65pc{
\delta\ar[d]^{\rotatebox{90}{$\sim$}}\\ \delta_s\ar[d]^H \\ \delta'_s\ar[d]^{\rotatebox{90}{$\sim$}} \\ \delta'
}\qquad
\xymatrix@=1.5pc{
a\ar[r]\ar[d]&b\ar[r]\ar[d]&c\ar[r]\ar[d]&a[1]\ar[d]\\
a\ar[r]\ar[d]^f&b\ar[r]\ar[d]^g&C(u)\ar[r]\ar[d]^h&a[1]\ar[d]^{f[1]}\\
a'\ar[r]\ar[d]&b'\ar[r]\ar[d]&C(u')\ar[r]\ar[d]&a'[1]\ar[d]\\
a'\ar[r]&b'\ar[r]&c'\ar[r]&a'[1]
}
\end{equation*}
where $H$ is of ``homotopy type'' which means in $K^*(\A)$ that if we fix representative maps
for the arrows then $h$ is upper triangular with the upper triangular entry being a homotopy
$a\to b'$. This can be extended to distinguished maps of $n$-triangles via the $2$-faces of the $n$-triangles. The verification of the axioms is somewhat tedious and the details can be found in the thesis \cite{Hus}. It should be stressed that the key idea
for this goes back to Neeman's 1991 paper in which a new set of axioms for triangulated categories is considered
which are better suited to the homotopy or derived category. In the introduction to that paper, for example, Neeman
shows that the Lightning Strike Axiom holds. The axioms are also likely to hold in similar situations where there
are standard models for distinguished triangles and maps between them (for example the stable homotopy category).
In the case where $\A$ is abelian we can deduce that $D^*(\A)$ is also fully triangulated via the canonical map $Q:K^*(\A)\to D^*(\A)$. It should not be too hard to generalise this to that stable homotopy category but we do not do that here. It is also possible that the axioms characterise these constructions and that these are essentially the only models. 
\end{example}
\begin{remark}\label{r:axioms}
  \begin{enumerate}
    \item In fact we can also include the case $n=1$. Then the axioms just say that $\ob(D_1)=\ob(\C)$ and
$\arr(D_1)=\arr{\C}$.
\item The Faces Axiom says that the face and degeneracy functors restrict to functors $D_n\to D_{n-1}$ and $D_{n-1}\to D_n$ respectively.
\item The final axiom is so named because Neeman called such maps ``lightning strikes'' and we
  will follow his convention here too. It will turn out that the converse of the Lightning
  Strike Axiom (\ref{a:ls}) also holds: if $(f,g,h)$ is distinguished then so it $(f,g,h+v'\tau
  w)$ for any map $\tau:\codom(w)\to\dom(v')$ but this needs some more machinery and we will see
  it in the next section. In the homotopy category, it is possible to formulate analogues of this axiom for $n>2$ but it seems likely that
  these would follow from the $n=2$ case and are messy to write down.
\item We do not require the induced map in the Bases Axiom
  (\ref{a:bases}) to be unique up to unique isomorphism. But the point is that these maps are
  considerably more special than a general map between the $n$-triangles. 
  \item The constructions in the Bases Axiom are not natural (still) but they are more constrained than the construction of maps
  of triangles or octahedra in a triangulated category. Consequently, the constructions are natural in certain contexts such as where
  the objects lie in a single abelian subcategory.
  \item It is natural to consider the class of all fully triangulated categories and give it the structure of a category where
  the functors are functors of categories preserving $D_n$. Perhaps these can be called
  \strong{fully exact functors} but we do not consider them further in this paper.
\end{enumerate}
\end{remark}
The following are easy consequences of the axioms.
\begin{prop}\label{p:axioms}
\begin{enumerate}
\item\label{isoms} Any object of $\Pi_n$ which is isomorphic to a distinguished $n$-triangle is again an distinguished $n$-triangle.
\item The Faces Axiom (\ref{a:faces}) holds for objects as well as morphisms. 
\item For any object $a$, we have $a=a\to0\to a[1]$ is a distinguished $2$-triangle.
\end{enumerate}
\end{prop}
 In particular it follows that a fully triangulated category is pre-triangulated in the usual
 sense. We shall see below that it is actually triangulated by Proposition \ref{p:strongdist} below.

\begin{lemma}\label{l:isom}
If $f$ in the Bases Axiom (\ref{a:bases}) is an isomorphism then any choice of map of associated
distinguished triangles restricting to $f$ on their bases is an isomorphism. 
\end{lemma}
\begin{proof}
By the Faces Axiom (\ref{a:faces}), every $2$-triangular face is distinguished and so the usual
argument with the five lemma and Yoneda's lemma for pre-triangulated categories implies that the
map on the   cones of $a_{0,i}\to a_{0,j}$ for $0<i<j\leq n$ are isomorphisms and hence so is
the extension of $f$ since every vertex of the $n$-triangles appears as a cone on a map of the
base. 
\end{proof}

It is well known that the Octahedral Axiom for triangulated categories is
equivalent, in a pre-triangulated category, to the existence of an octahedral diagram on a base
of two composable maps. In fact, this works for all $n$: 
\begin{prop}\label{p:strongdist}
In any fully triangulated category, given an $n$-dimensional simplicial diagram $\sigma$ in $\C$
and choices of distinguished $2$-triangles $\delta_{ij}$: $a_i\to a_j\to a_{i,j}\to a_i[1]$ on
each map $a_i\to a_j$ of $\sigma$, there is a distinguished $n$-triangle $\delta$ with base
$\sigma$ and containing each $\delta_{ij}$ as faces. 
\end{prop} 

\begin{proof}
Observe that the data in the hypothesis includes all of the vertices of $\delta$ and some of the
edges, namely $e_{0,i,j}$ and $e_{i,0,j}$ for $0<i<j\leq n$. We prove it by induction on
$n$. The $n=2$ case is obvious from the Bases Axiom (5). We assume that the proposition holds
for $n-1$. Let $\delta$ be a distinguished $n$-triangles on the given base and consider its face
$F_n$. Then by previous Lemma \ref{l:isom} this is isomorphic to the distinguished
$n-1$-triangle $\delta'$ which exists by induction formed from the subset of triangles of the
hypothesis with $j<n$. Then we can extend this isomorphism to an isomorphism of the $n$-th base
$\sigma_n$ of $\delta$ by choosing automorphisms $\phi_i$ of $a_{i,n}$ for $1\leq i<n$ which
extend the identity map $a_{0,i}\to a_{0,n}$ to a distinguished map of triangles. Now we can
define the missing maps $a_{i,n-1}\to a_{i,n}$ and $a_{i,n}\to a_{i+1,n}$ by composition with
these automorphisms. It is an easy exercise to check that this gives an isomorphism of rectified
simplex diagrams and since one is distinguished, so is the other by Proposition
\ref{p:axioms}(\ref{isoms}). 
\end{proof}
As illustration of this is given in the proof of Proposition \ref{p:facemaps} below.
\begin{remark}
In particular, all of the usual consequences of the axioms of triangulated categories follow. For
example, there is a canonical additive structure determined by the triangulated structure (the biproduct
$a\oplus b$ is a cone on the zero map $a[-1]\to b$ and $0$ is the cone on any identity map).  
\end{remark}
As an application of this proposition we show that the maximal cycle of maps of faces are
distinguished maps of $n-1$-triangles. 
\begin{prop}\label{p:facemaps}
Suppose $\delta$ is a distinguished $n$-triangle. Then each maps of faces $F_n\to
F_{n-1}\to\cdots\to F_0\to \tau^n F_n$ is distinguished. 
\end{prop}
\begin{proof}
We describe the $3$-triangle case in detail first. Suppose we are given a $3$-triangle
\[\xymatrix@=1.4pc{
a\ar[r]^u&b\ar[r]^f\ar[d]^v&c\ar[d]^{v'}\\
&d\ar[r]^g\ar@{~>}[ul] &e\ar[d]^{g'}\ar[r]^{w'}&a[1]\ar[d]^{u[1]}\\
&&p\ar[r]^{f''}\ar@{~>}[ul]&b[1]
}\]
Then we have a distinguished octahedron by Proposition \ref{p:strongdist}
\[\xymatrix@=2.5pc{a\ar[r]^{fu} & c\ar@{=}[r]\ar[d]^{v'}& c\ar[d]^{v'}\\
&e\ar[r]^{\beta}\ar@{~>}[ul] &e\ar[r]^{w'}\ar[d]^{0}&a[1]\ar[d]^{fu[1]}\\
&&0\ar[r]^{0}\ar@{~>}[ul]&c[1]
}\]
For some map $\beta$ as indicated. Note that $\beta$ is an isomorphism as $e\sra{\beta}e\to0\to
e[1]$ is distinguished. We also have the horizontal composite is $w'$ and so $w'\beta=w'$. The
middle square gives $\beta v'=v'$.  

Now consider the map of bases:
\[\xymatrix@=1.8pc{
a\ar[r]^u\ar@{=}[d] &b\ar[r]^f\ar[d]^{f} &c\ar@{=}[d]\\
a\ar[r]^{fu} &c \ar@{=}[r]& c
}\]
By the Bases Axiom this lifts to a map of distinguished octahedra via maps $\zeta:d\to e$,
$\alpha:e\to e$ and $0:p\to 0$. We now consider the restriction to the faces in turn: 
\[\xymatrix@=1.9pc{
d\ar[r]^g\ar[d]_\zeta &e\ar[r]^{g'}\ar[d]^<<<<{\alpha}& p\ar[r]^{v[1]f''}\ar[d]^\eta&d[1]\ar[d]^{\zeta[1]}\\
e\ar[r]^{\beta}&e\ar[r]&0\ar[r]&e[1]
}\]
From this we have $\zeta[1] v[1]f''=0$ and $\alpha\zeta=\beta g$. We also have
\[\xymatrix@=1.9pc{
a\ar@{=}[d]\ar[r]^{fu} & c\ar[r]^{v'}\ar@{=}[d]& e\ar[r]^{w'}\ar[d]^{\alpha}&a[1]\ar@{=}[d]\\
a\ar[r]^{fu}&c\ar[r]_{v'} & e\ar[r]^{w'} &a[1]
}\]
From this we have $w'\alpha=w'$ and $\alpha v'=v'$. Then $(1,1,\alpha)$ is a map of
distinguished 2-triangles and so $\alpha$ is an isomorphism. 
Finally, we have
\[\xymatrix@=1.9pc{
a\ar[r]^u\ar@{=}[d] & b\ar[r]^v\ar[d]^f &d\ar[r]^{w'g}\ar[d]^{\zeta} &a[1]\ar@{=}[d]\\
a\ar@{=}[d]\ar[r]^{fu} &c\ar@{=}[d]\ar[r]^{v'} &e\ar[d]^{\beta^{-1}\alpha}\ar[r]^{w'}&a[1]\ar@{=}[d]\\
a\ar[r]^{fu} &c\ar[r]^{v'} &e\ar[r]^{w'}&a[1]
}\]
The top map of triangles is distinguished by the Faces Axiom and the diagram commutes via the
identities we have found. The bottom map of triangles is an isomorphism and so by the
Isomorphism Axiom this composite is a distinguished map. But the composite is exactly the map of
the faces since $g=\beta^{-1}\alpha\zeta$. 

It follows that the other face maps are distinguished by rotational symmetry.

The general case follows similarly by observing that the map $F_n\to F_{n-1}$ is the identity on
$F_n\cap F_{n-1}$ and $n-1$ maps $a_{i,n-1}\to a_{i,n}$ for $0\leq i<n-1$. 
Given our standard $n$-triangle with vertices $a_{i,j}$ for $0\leq i<j\leq n$. We consider
another distinguished $n$-triangle  
\begin{equation}\label{eq:degex}
\xymatrix@=1.5pc{
\cdots&a_{0,n-2}\ar[r]& a_{0,n}\ar@{=}[r]\ar[d]^{f_0} & a_{0,n}\ar[d]^{f_0}\\
&\cdots&a_{1,n}\ar[r]^{\beta_1}\ar[d]^{f_1}&a_{1,n}\ar[r]\ar[d] &a_{0,1}[1]\ar[d]\\
&\ddots&\vdots\ar[d]^{f_{n-2}}&\vdots\ar[d]&\vdots\ar[d]\\
&&a_{n-3,n}\ar[r]^{\beta_{n-3}}\ar[d]^{f_{n-3}} &a_{n-3,n}\ar[r]\ar[d] & a_{0,n-3}[1]\ar[d]\\
&&a_{n-2,n}\ar[r]^{\beta_{n-2}}\ar@{~>}[ul]&a_{n-2,n}\ar[r]\ar[d] &a_{0,n-2}[1]\ar[d]\\
&&&0\ar[r] \ar@{~>}[ul]&a_{0,n}[1]
}\end{equation}
Since the 2-triangle faces are distinguished we see that the $\beta_i$ maps are all
isomorphisms. Note also that the $n^\text{th}$ face of this $n$-triangle is $F_{n-1}$  
Now consider the map of bases:
\[\xymatrix@=1.5pc{
a_{0,1}\ar@{=}[d]\ar[r]&a_{0,2}\ar@{=}[d]\ar[r]&\cdots\ar[r] &a_{0,n-2}\ar[rr]^{u_{n-2}}\ar@{=}[d]&&a_{0,n-1}\ar[r]^{u_{n-1}}\ar[d]_{u_{n-1}}&a_{0,n}\ar@{=}[d]\\
a_{0,1}\ar[r] &a_{0,2}\ar[r]&\cdots\ar[r] &a_{0,n-2}\ar[rr]_{u_{n-1}u_{n-2}} &&a_{0,n}\ar@{=}[r] &a_{0,n}
}\]
Then we have a distinguished map $\phi:F_n\to F_{n-1}$ by restriction. As for the $n=3$ case we
then have maps $\alpha_i:a_{i,n}\to a_{i,n}$ for $1\leq i\leq n-2$ and again these must be
isomorphisms. Then we have that the original map of faces factors through $\phi$ and an
automorphism $F_{n-1}$. The other face maps are also distinguished by rotational symmetry
again. 
\end{proof}
In fact the $n$-triangle  (\ref{eq:degex}) is actually a degeneration of the $F_{n-1}$ face in
which the $\beta_i$ are the identity maps but we do not need this. 

There is a weak converse to this proposition.
\begin{prop}\label{p:mapfaces}
Given a distinguished map between two distinguished $n-1$-triangles $A$ and $B$ sharing a common
$n-2$-triangle and such that the map restricts to the identity on the intersection, there is an
$n$-triangle $\delta$ with faces $F_n=A$ and $F_{n-1}=B$.   
\end{prop}
The proof is complicated and we omit the details. Instead we illustrate the technique with the $n=3$ case.
\begin{prop} Given a distinguished map 
\[\xymatrix@=1.8pc{a\ar[r]^u\ar@{=}[d] &b\ar[r]^v\ar[d]^g &c\ar[r]^w \ar[d]^h & a[1]\ar@{=}[d]\\
a\ar[r]^{gu} &b'\ar[r]^{v'} &c'\ar[r]^{w'} &a[1]
}\]
of distinguished $2$-triangles, there is a distinguished $3$-triangle containing this as a map of two faces.
\end{prop}
\begin{proof}
By Proposition \ref{p:strongdist} there is a distinguished $3$-triangle
\begin{equation}\label{eq:face}
\xymatrix@=1.8pc{
a\ar[r]^u & b\ar[r]^g\ar[d]^v &b'\ar[d]^{v'}\\
& c\ar[r]^{h'} \ar@{~>}[ul]^w& c'\ar[r]^{w'}\ar[d]^{q} &a[1]\ar[d]^{u[1]}\\
&& c''\ar[r]^\psi \ar@{~>}[ul]^{v[1]\psi}& b[1]
}\end{equation}
By Proposition \ref{p:facemaps}, $(1,g,h')$ is also distinguished and so there is a map $\theta:a[1]\to b'$ such that $h'-h=v'\theta w=v'\theta w'h$. So $h'=\alpha h$, where $\alpha=1+v'\theta w'$. Then $w'\alpha=w'$ and $\alpha v'=v'$ and so $(1,1,\alpha)$ is an isomorphism and hence so is $\alpha$. Then we have an isomorphism of 2-triangles
\[\xymatrix@=1.8pc{
c\ar@{=}[d] \ar[r]^h &c'\ar[d]^\alpha \ar[r]^{q\alpha} & c'' \ar@{=}[d] \ar[r]^{v[1]\psi} & c[1]\ar@{=}[d]\\
c\ar[r]^{h'} &c'\ar[r]^q &c''\ar[r]^{v[1]\psi} &c[1]
}\]
Then the top triangle is distinguished. Moreover, we have a 3-triangle
\[\xymatrix@=1.8pc{
a\ar[r]^u & b\ar[r]^g\ar[d]^v &b'\ar[d]^{v'}\\
& c\ar[r]^{h} \ar@{~>}[ul]^w& c'\ar[r]^{w'}\ar[d]^{q\alpha} &a[1]\ar[d]^{u[1]}\\
&& c''\ar[r]^\psi \ar@{~>}[ul]^{v[1]\psi}& b[1]
}\]
This is isomorphic to (\ref{eq:face}) via the identity everwhere except $c'$ where we have $\alpha$ and so is also distinguished.
\end{proof}

\section{Consequences of the axioms}
We shall now consider some key consequences of the axioms.

\subsection{Direct sums}
Recall that in a pre-triangulated category two $2$-triangles are distinguished if and only if
their direct sum is distinguished 
(see, for example, \cite[Propositions 1.2.1 \&\ 1.2.3]{NeemanBook}). We now show that this also
works for direct sums of distinguished maps. We only consider the $n=2$ case. Our set up will
be two maps of $2$-triangles: 
\begin{equation}\label{eq:triangpair}
\vcenter{\vbox{\hbox{\xymatrix@=1.6pc{
a_1\ar[r]^{u_1}\ar[d]^{f_1} &b_1\ar[r]^{v_1}\ar[d]^{g_1} &c_1\ar[r]^{w_1}\ar[d]^{h_1} &a_1[1]\ar[d]^{f_1[1]}\\
a'_1\ar[r]^{u'_1} &b'_1\ar[r]^{v'_1} &c'_1\ar[r]^{w'_1} &a'_1[1]
}\qquad
\xymatrix@=1.6pc{
  a_2\ar[r]^{u_2}\ar[d]^{f_2} &b_2\ar[r]^{v_2}\ar[d]^{g_2} &c_2\ar[r]^{w_2}\ar[d]^{h_2} &a_2[1]\ar[d]^{f_2[1]}\\
a'_2\ar[r]^{u'_2} &b'_2\ar[r]^{v'_2} &c'_2\ar[r]^{w'_2} &a'_2[1]
 }}} }\end{equation}
and their direct sum:
\[
\xymatrix@=1.9pc{
a_1\oplus a_2\ar[r]^{u_1\oplus u_2}\ar[d]^{f_1\oplus f_2} &b_1\oplus b_2\ar[r]^{v_1\oplus
  v_2}\ar[d]^{g_1\oplus g_2} &c_1\oplus c_2\ar[r]^{w_1\oplus w_2}\ar[d]^{h_1\oplus h_2}
&a_1[1]\ar[d]^{(f_1\oplus f_2)[1]}\\
a'_1\oplus a'_2\ar[r]^{u'_1\oplus u'_2} &b'_1\oplus b'_2\ar[r]^{v'_1\oplus v'_2} &c'_1\oplus
c'_2\ar[r]^(0.6){w'_1\oplus w'_2} & **[r](a'_1\oplus a'_2)[1]
}\]
We shall abberviate these to maps $G_1$ and $G_2$ with direct sum $G_1\oplus G_2$.  
We aim to prove
\begin{thm}\label{t:sum}
The maps $G_1$ and $G_2$ are distinguished if and only if $G_1\oplus G_2$ is distinguished.
\end{thm}
From this we conclude that the sum of two distinguished maps is again distinguished:
\begin{cor}
If $G_1$ and $G_2$ are a distinguished map of the same two $2$-triangles then their sum $G_1+G_2$ is also distinguished. 
\end{cor}
\begin{proof}
This follows immediately by writing $G_1+G_2=\nabla(G_1\oplus G_2)\Delta$ via the canonical diagonal inclusion and projection maps $\Delta$ and $\nabla$.
\end{proof}
This contrast with the case of good maps.
\begin{prop}\label{p:sum}
For any $f_1$, $f_2$ ,$g_1$ and $g_2$, these can be completed to distinguished maps $G_1$ and $G_2$ such that  $G_1\oplus G_2$ is distinguished.
\end{prop}
The key to the proof is a distinguished $5$-triangle.
\begin{lemma}\label{l:5trisum}
  For any two distinguished 2-triangles $a_i\sra{u_i}b_i\sra{v_i}c_i\sra{w_i}a_i[1]$ for
  $i=1,2$, the following is a distinguished 5-triangle 
 \begin{equation}\label{eq:5tri}
  \vcenter{\xymatrix@=1.8pc{
    a_2\ar[r]^{\smat{0\\1}}& a_1\oplus a_2\ar[r]^{u_1\oplus 1}\ar[d]^{(1\;0)} &
                   b_1\oplus a_2 \ar[r]^{1\oplus u_2}\ar[d]^{(1\;0)} & b_1\oplus
                   b_2\ar[r]^{(0\;1)}\ar[d]^{1\oplus v_2} & b_2\ar[d]^{v_2} \\
  & a_1\ar[r]^{u_1}\ar@{~>}[ul]^0  & b_1\ar[r]^{\smat{1\\0}} \ar[d]^{v_1} &  b_1\oplus
                   c_2\ar[r]^{(0\;1)}\ar[d]^{v_1\oplus 1} &    c_2\ar[r]^{w_2}\ar[d]^{\smat{0\\1}} &  a_2[1]\ar[d]^{\smat{0\\1}}\\
  &&  c_1\ar[r]^{\smat{1\\0}}\ar@{~>}[ul]^{w_1} & c_1\oplus c_2\ar[r]^{w_1\oplus 1}\ar[d]^{(0\;1)} & a_1[1]\oplus
                        c_2\ar[r]^{1\oplus w_2}\ar[d]^{u_1\oplus 1}  & (a_1\oplus
                        a_2)[1]\ar[d]^{u_1[1]\oplus 1} \\
      &&&c_2\ar[r]^{\smat{0\\1}}\ar@{~>}[ul]^0 & b_1[1]\oplus c_2\ar[r]^{1\oplus w_2}\ar[d]^{(1\;0)} &
                        (b_1\oplus a_2)[1]\ar[d]^{1\oplus u_2[1]}\\
  &&&& b_1[1]\ar[r]^{\smat{1\\0}}\ar@{~>}[ul]^0 & (b_1\oplus b_2)[1]
                                               }}
\end{equation}
\end{lemma}
\begin{proof}
  We use Proposition \ref{p:strongdist} using ten 2-triangles which are easily seen to be distinguished:
  \begin{align}
    a_2&\lra{\smat{0\\1}} a_1\oplus a_2\lra{(1\;0)}a_1\lra{0}a_2[1]\\
    a_2&\lra{\smat{0\\1}}b_1\oplus a_2\lra{(1\;0)} b_1\lra{0}b_1[1]\\
    a_2&\lra{\smat{0\\u_2}} b_1\oplus b_2\lra{1\oplus v_2} b_1\oplus c_2\lra{(0\; w_2)} a_2[1]\\
    a_2&\lra{\smat{0\\u_2}} b_2\lra{v_2} c_2\lra{w_2}a_2[1] \label{eq:second}\end{align}
\begin{align}
    a_1&\oplus a_2\lra{u_1\oplus 1}b_1\oplus a_2\lra{(v_1\;0)} c_1\lra{\smat{w_1\\0}}(a_1\oplus
    a_2)[1]\\
    a_1&\oplus a_2\lra{u_1\oplus u_2}b_1\oplus b_2\lra{v_1\oplus v_2} c_1\oplus c_2\lra{w_1\oplus
      w_2} (a_1\oplus a_2)[1]\label{eq:sum} \end{align}
 \begin{align}     
    a_1&\oplus a_2\lra{(0\;u_2)}b_2\lra{\smat{0\\v_2}} a_1[1]\oplus c_2\lra{1\oplus w_2}
    (a_1\oplus a_2)[1]\\
    b_1&\oplus a_2\lra{1\oplus u_2}b_1\oplus
    b_2\lra{(0\;v_2)}c_2\lra{\smat{0\\w_2}}(b_1\oplus a_2)[1]\\
    b_1&\oplus a_2\lra{(0\;u_2)} b_2\lra{\smat{0\\v_2}}b_1[1]\oplus c_2\lra{1\oplus
      w_2}(b_1\oplus a_2)[1]\\
    b_1&\oplus b_2\lra{(0\;1)} b_2\lra{0}b_1[1]\lra{\smat{1\\0}}(b_1\oplus b_2)[1]
  \end{align}
  This results in a distinguished 5-triangle:
   \begin{equation}
  \vcenter{\xymatrix@=1.8pc{
    a_2\ar[r]^{\smat{0\\1}}& a_1\oplus a_2\ar[r]^{u_1\oplus 1}\ar[d]^{(1\;0)} &
                   b_1\oplus a_2 \ar[r]^{1\oplus u_2}\ar[d]^{(1\;0)} & b_1\oplus
                   b_2\ar[r]^{(0\;1)}\ar[d]^{1\oplus v_2} & b_2\ar[d]^{v_2} \\
  & a_1\ar[r]^{A}\ar@{~>}[ul]^0  & b_1\ar[r]^{B} \ar[d]^{E} &  b_1\oplus
                   c_2\ar[r]^{C}\ar[d]^{H} &    c_2\ar[r]^{w_2}\ar[d]^{I} &  a_2[1]\ar[d]^{\smat{0\\1}}\\
  &&  c_1\ar[r]^{J}\ar@{~>}[ul]^{w_1} & c_1\oplus c_2\ar[r]^{K}\ar[d]^{L} & a_1[1]\oplus
                        c_2\ar[r]^{1\oplus w_2}\ar[d]^{M}  & (a_1\oplus
                        a_2)[1]\ar[d]^{u_1[1]\oplus 1} \\
      &&&c_2\ar[r]^{N}\ar@{~>}[ul]^0 & b_1[1]\oplus c_2\ar[r]^{1\oplus w_2}\ar[d]^{P} &
                        (b_1\oplus a_2)[1]\ar[d]^{1\oplus u_2[1]}\\
  &&&& b_1[1]\ar[r]^{\smat{1\\0}}\ar@{~>}[ul]^0 & (b_1\oplus b_2)[1]
                                               }}
\end{equation}
Then $A=u_1$, $B=\smat{1\\ 0}$, $C=(0\;\alpha)$, for some $\alpha$. Then $\alpha v_2=v_2$,
$w_2\alpha=w_2$ and so $\alpha$ is an isomorphism. So we can apply its inverse at the codomain
of $C$ and identity everywhere else to give a new 5-triangle (if necessary changing $I$) and we
can assume $\alpha=1$. Now we
see that $E=v_1\beta$ for some isomorphism $\beta$, and  $I=\smat{0\\ \eta}$,
where $w_2\eta=w_2$, and $\eta v_2=v_2$ so that $\eta$ is also an isomorphism.

Now $N=\smat{\gamma\\ \theta}$ and $P=(1\;\kappa)$. We have $\theta$ is an ismorphism which we
can also take to be the identity (up to an isomorphism of the 5-triangle fixing the other maps) and
then $\gamma=-\kappa$. But $(1\oplus w_2)N=\smat{0\\w_2}$ and so $\gamma=-\kappa=0$.
Up to isomorphism (possibly changing $H$, $K$ and $L$) we can assume $J=\smat{1\\ 0}$ and $\beta=1$. 

Now let $H=\smat{x&y\\m&n}$, $K=\smat{p&q\\t&z}$, $L=(0\;\zeta)$  and
$M=\smat{r&s\\\pi&\xi}$. Here $\zeta$ is an isomorphism. Then
$HB=\smat{x\\m}=\smat{v_1\\0}$ and $H(1\oplus v_2)=\smat{v_1&yv_2\\0&nv_2}=\smat{v_1&0\\0&v_2}$
and so $yv_2=0$, $nv_2=v_2$ and $x=v_1$. Similarly, from $K(1\oplus w_2)=w_1\oplus w_2$ we have $q=0$,
$p=w_1$, $w_2t=0$ and $w_2z=w_2$. Now $KH=\smat{0&w_1y\\ty_1&ty+zn}=IC=\smat{0&0\\0&\eta}$ and so
$w_1y=0$, $tv_1=0$ and $ty+zn=\eta$. We have $(1\oplus w_2)M=u_1[1]\oplus w_2$ and so $r=u_1[1]$,
$s=0$, $w_2\pi=0$ and $w_2\xi=w_2$. From $MIv_2$ we have $\xi v_2=v_2$. Then $\xi$ is an
isomorphism. Then from $MK$ we have $z=1$ and
$t=-\xi^{-1}\pi w_1$. From $(1\oplus w_2)NLH=(u_1[1]\oplus1)\smat{0\\1}w_2C$ we have
$w_2\zeta n=w_2$. From $LH(1\oplus v_2)$ we have $\zeta n v_2=v_2$.  So $n$
  is also an isomorphism. We can now construct an
  isomorphism of the 5-triangle with isomorphisms $\zeta^{-1}$, $\smat{1&-y\\0&n^{-1}}$ and
  $\smat{1&0\\-\pi&\xi^{-1}}$ at $c_2$ (in position $a_{3,4}$), $c_1\oplus c_2$ and $a_1[1]\oplus c_2$ respectively. Then, since
  $ty+n=\eta$ we have an isomorphism to our given 5-triangle of the lemma which must therefore be distinguished. 
\end{proof}
We can now prove the proposition.
\begin{proof}
We have two 5-triangles, one for the domains of $G_1$ and $G_2$ which we will call $\delta$
and is exactly the one in (\ref{eq:5tri}) and the other for the
codomains (the ``dashed'' triangles) which we shall denote $\delta'$. Then the key observation is that the original 2-triangles
appear in the 5-triangle as the faces (\ref{eq:second}), (\ref{eq:sum}) and $a_{1,2}\to
a_{1,3}\to a_{2,3}\to a_{1,2}[1]$. A map between the bases of $\delta$ and $\delta'$ is exactly
given by the maps $f_1$, $g_1$, $f_2$ and $g_2$ of (\ref{eq:triangpair}). Then if we extend this
to a distinguished map $\Xi:\delta\to\delta'$ this will restrict to distinguished maps of all faces
including our three 2-triangles. Let the maps at $c_1$ and $c_2$ be $h_1$ and $h_2$
respectively. Now suppose the map at $c_1\oplus c_2$ is given by
$A$. Then $A\smat{1\\0}=\smat{h_1\\0}$, $(0\;1)A=(0\;h_2)$ and
$A(v_1\oplus 1)=\smat{v'_1h_1&0\\0&h_2}$ and so $A=h_1\oplus h_2$.  This implies that there is
at least one choice of distinguished maps $(f_i,g_i,h_i)$ such that $(f_1\oplus f_2,g_1\oplus
g_2,h_1\oplus h_2)$ is distinguished.
\end{proof}
 To complete the proof of Theorem \ref{t:sum} we need to understand better the class of distinguished
 maps $G_i$ for fixed $f_1$, $f_2$, $g_2$ and $g_2$.
 
\subsection{Good maps, lightning strikes and the Neeman octahedron}
So far we have not used the Lightning Strike Axiom (\ref{a:ls}) and this will be the key to understanding the problem above.
\begin{dfn}
Given two $2$-triangles, $a\sra{u} b\sra{v} c\sra{w} a[1]$ and $a'\sra{u'} b'\sra{v'}c'\sra{w'}
a'[1]$, we define a relation on maps  $(f,g,h)$ between 
them by $(f,g,h)\sim (f,g,h')$ if $h'-h=v'\tau w$ for some $\tau:a[1]\to b'$. In other words, $h$
and $h'$ differ by a lightning strike. It is an easy 
exercise to check that this is an equivalence relation and we call it the \strong{lightning
  strike equivalence relation}.  
\end{dfn}
Our Lightning Strike Axiom tells us that any two distinguished maps $(f,g,h)$ and $(f,g,h')$
are lightning strike equivalent. It will turn out that the lightning strike equivalence class
containing a distinguished map consists only of distinguished maps.

Neeman also considered lightning strikes in the context of his so called ``good'' maps (see
\cite{Neeman1991}.
Recall that a map of triangles $(f,g,h):(a\sra{u} b\sra{v}c\sra{w} a[1]) \lra{} (a'\sra{u'}
b'\sra{v'}c'\sra{w'} a'[1])$ is \strong{good} if the mapping cone
\[\xymatrix@=3pc{
a'\oplus b\ar[r]^{\smat{u'&g\\0&-v}} & b'\oplus c\ar[r]^(0.45){\smat{v'&h\\0&-w}}& c'\oplus
a[1]\ar[r]^(0.45){\smat{w'&f[1]\\0&-u[1]}}& (a'\oplus b)[1]  } 
\]
is distinguished. The notion of good is clearly preserved by isomorphism of 2-triangles.
He also showed that if $(f,g,h)$ is good and $\tau:a[1]\to b'$ then $(f,g,h+w'\tau v)$ is also
good. But good maps do not necessarily form a single lightning strike equivalence class and
they also do not necessarily compose to give good maps.
He also showed that any commuting square can be completed to a good map of distinguished
triangles. To do so he used a special octahedron which we will see below and we will call the
\strong{Neeman octahedron}. He also showed that a good map of the form $(1,g,h)$ can be
completed to an octahedron. This is also true of distinguished maps by Proposition
\ref{p:mapfaces}. But for good maps it may not be true that the other maps of faces are good.

In this section we aim to use Neeman's techniques to complete the proof of Theorem \ref{t:sum},
show that the converse of the Lightning Strike Axiom holds (so that the collection of
distinguished maps of 2-triangles with first two maps fixed forms a single lightning strike equivalence class),
and that distinguished maps of 2-triangles are always good.

\begin{dfn}
For a map of 2-triangles \[G=(f,g,h):(a\sra{u} b\sra{v}c\sra{w} a[1]) \lra{} (a'\sra{u'}
b'\sra{v'}c'\sra{w'} a'[1]),\] we consider the following map of triangles:
\[\xymatrix@=2.2pc{
a'\oplus a\ar[r]^(0.45){\smat{1&0\\0&u\\0&0}}\ar@{=}[d] & a'\oplus b\oplus
b'\ar[r]^(0.6){\smat{0&0&1\\0&v&0}}\ar[d]^{(u'\;g\;1)} & b'\oplus
c\ar[r]^(0.4){\smat{0&0\\0&w}}\ar[d]^(0.4){\smat{v'&h\\0&-w}} &(a'\oplus a)[1]\ar@{=}[d]\\
a'\oplus a\ar[r]^{(u'\;gu)} & b'\ar[r]^{\smat{v'\\0}} & c'\oplus a[1]\ar[r]^{\smat{w'&f\\0&-1}} &
  (a'\oplus a)[1]
  }\]
We shall call it \strong{Neeman's map} (see the proof of Theorem 1.8 in \cite{Neeman1991}) and
denote it $N(G)$.
\end{dfn}
What this does is to convert a general map of triangles $G$ into a map of triangles in which the
first map is the identity (and so is a candidate for a map of faces of an octahedron). Neeman
proves that if $N(G)$ can be completed to an octahedron then  $G$ is good. The key to this
section is that $N(G)$ is distinguished if and only of $G$ is distinguished. 
\begin{lemma}
  $N(G)$ is isomorphic to $N'(G)$ defined to be
  \[\xymatrix@=2.2pc{
a'\oplus a\ar[r]^(0.45){\smat{1&0\\0&u\\0&0}}\ar[d]^{\smat{1&f\\0&-1}} & a'\oplus b\oplus
b'\ar[r]^(0.6){\smat{0&0&1\\0&v&0}}\ar[d]^{(u'\;g\;1)} & b'\oplus
c\ar[r]^(0.4){\smat{0&0\\0&w}}\ar[d]^(0.4){\smat{v'&h\\0&-w}} &(a'\oplus a)[1]\ar[d]^{\smat{1&f[1]\\0&-1}}\\
a'\oplus a\ar[r]^{(u'\;0)} & b'\ar[r]^{\smat{v'\\0}} & c'\oplus a[1]\ar[r]^{\smat{w'&0\\0&1}} &
  (a'\oplus a)[1]
  }\]
 \end{lemma}
\begin{proof}
The isomorphism is given by $\displaystyle\left(\mat{1&f\\0&-1},1,1\right)$.
\end{proof}
So $N(G)$ is distinguished if and only if so is $N'(G)$. We will also consider another related
map $N''(G)$ of distinguished triangles:
  \[\xymatrix@=2.2pc{
a'\oplus a\ar[r]^(0.45){\smat{1&0\\0&u\\0&0}}\ar[d]^{(1\;f)} & a'\oplus b\oplus
b'\ar[r]^(0.6){\smat{0&0&1\\0&v&0}}\ar[d]^{(u'\;g\;1)} & b'\oplus
c\ar[r]^(0.4){\smat{0&0\\0&w}}\ar[d]^{(v'\;h)} &(a'\oplus a)[1]\ar[d]^{(1\;f[1])}\\
a'\ar[r]^{u'} & b'\ar[r]^{v'} & c'\ar[r]^{w'} &
  a'[1]
  }\]
  All three of these share a useful property:
  \begin{lemma}
Fix the first two maps of $N(G)$, $N'(G)$ and $N''(G)$. Then any distinguished completions
take the form $\mat{v'&h\\0&-w}$ for $N(G)$ and $N'(G)$, and $(v',h+v'\tau w)$ for $N''(G)$ for some $\tau:a[1]\to b'$ and fixed map $h$. Moreover, any two distinguished completions are isomorphic via an automorphism of the top triangles.
  \end{lemma}
  \begin{proof}
The first statement is clear by considering a general lightning strike in each case. For the
second, we observe that $\left(1,1,\mat{1&\tau w\\0&1}\right)$ is the required isomorphism.
  \end{proof}
Since distinguished maps are closed under composition with isomorphisms we see that every
element of the lightning equivalence class containing one distinguished map is
distinguished for $N(G)$, $N'(G)$ and $N''(G)$. We now show this also holds for $G$ itself.
\begin{prop}\label{p:Nddash}
$N''(G)$ is distinguished if and only if $G$ is distinguished.
\end{prop}
\begin{proof}
By the preceding observation, it suffices to show that there is at least one map $h$ for which
both $G$ and $N''(G)$ are distinguished. We proceed as we did for Proposition \ref{p:sum} by
  showing that both $G$ and $N''(G)$ are restrictions of a single map of distinguished
  triangles. The distinguished triangles we need are 5-triangles again. The domain 5-triangle is
  \[\xymatrix@=3pc{
    a \ar[r]^(0.4){\smat{0\\1}} & a'\oplus a\ar[r]^(0.4){\smat{1&0\\0&u\\0&0}} \ar[d]^(0.45){(1\;0)} & a'\oplus
    b\oplus b'\ar[r]^(0.55){\smat{0&1&0\\u'&0&1}} \ar[d]^(0.45){\smat{1&0&0\\0&v&0\\u'&0&1}} & b\oplus b'
    \ar[r]^{(1\;0)} \ar[d]^{v\oplus 1} & b\ar[d]^{v} \\
    & a' \ar@{~>}[ul]^0 \ar[r]^(0.35){\smat{1\\0\\u'}} & a'\oplus c\oplus b' \ar[r]^(0.55){\smat{0&1&0\\0&0&1}}
    \ar[d]^(0.45){\smat{-u'&0&1\\0&1&0}} & c\oplus b' \ar[r]^{(1\;0)} \ar[d]^{1\oplus v'} &c \ar[r]^{w} \ar[d]^{\smat{0\\1}}
    & a[1]\ar[d]^{\smat{0\\1}} \\
    && b'\oplus c\ar@{~>}[ul]^0 \ar[r]^{\smat{0&1\\v'&0}}& c\oplus
    c'\ar[r]^(0.45){\smat{0&w'\\1&0}}\ar[d]^{(0\;w')} &a'[1]\oplus
    c\ar[d]^{\smat{1&0\\u'&0}}\ar[r]^{1\oplus w}
    &(a'\oplus a)[1]\ar[d]\\
    &&& a'[1]\ar@{~>}[ul]^{\smat{-u'\\0}} \ar[r]^(0.4){\smat{1\\-u'}} & (a'\oplus
    b')[1]\ar[r]^(0.45){\smat{1&0\\0&0\\0&1}} \ar[d]^{(u'\;1)}& (a'\oplus  b\oplus b')[1] \ar[d]\\ 
    &&&& b'[1]\ar@{~>}[ul]^0\ar[r]^{\smat{0\\1}} &(b\oplus b')[1]    
  }\]
  and the codomain is the following degenerate 5-triangle 
  \[\xymatrix@=2pc{
a'\ar@{=}[r] & a'\ar[r]^{u'}\ar[d] & b'\ar@{=}[r] \ar[d]^{v'} & b'\ar@{=}[r]\ar[d]^{v'} &
b'\ar[d]^{v'}\\
& 0\ar@{~>}[ul]\ar[r] &c'\ar@{=}[r]\ar@{=}[d] &c'\ar@{=}[r] \ar@{=}[d] &c'\ar[r]^{w'}\ar@{=}[d]
&a'[1]\ar@{=}[d]\\
&& c'\ar@{=}[r] \ar@{~>}[ul] & c'\ar@{=}[r]\ar[d] & c'\ar[d] \ar[r]^{w'} &a'[1]\ar[d]^{u'[1]}\\
&&&0\ar@{~>}[ul] \ar@{=}[r] & 0\ar[r]\ar@{=}[d] & b'[1]\ar@{=}[d] \\
&&&& 0\ar[r]\ar@{~>}[ul] &b'[1]
  }\]
The first of these is distinguished using the same method as in the proof of Lemma
\ref{l:5trisum}. Then we consider the following map of bases:
\[\xymatrix@=2pc{
   a \ar[r]^(0.4){\smat{0\\1}}\ar[d]^f & a'\oplus a\ar[r]^(0.4){\smat{1&0\\0&u\\0&0}} \ar[d]^{(1\;f)} & a'\oplus
    b\oplus b'\ar[r]^(0.55){\smat{0&1&0\\u'&0&1}} \ar[d]^{(u'\;g\;1)} & b\oplus b'
    \ar[r]^(0.6){(1\;0)} \ar[d]^{g\oplus 1} & b\ar[d]^{g} \\
 a'\ar@{=}[r] & a'\ar[r]^{u'}& b'\ar@{=}[r] & b'\ar@{=}[r]&
b'   }\]
Then pick a distinguished map extending this base map. Then the map at $b'\oplus c$ take sthe
form $\mat{v'&h\\0&-w}$ and then the map at $a_{1,5}$ is also $h$. But then we have $(f,g,h)$ as
a map of the face $a_{0,1}\to a_{0,5}\to a_{1,5}\to a_{0,1}[1]$ which is $G$. So we have, for
    this choice of $h$, that both $G$ and $N''(G)$ are distinguished as required.
\end{proof}
Since the lightning equivalence classes now coincide we see that $G$ is distinguished if and
only if $N''(G)$ is distinguished.
This also proves the converse of Axiom (6):
\begin{prop}
If $(f,g,h)$ is distinguished then so it $(f,g,h+v'\tau w)$ for any map $\tau:\codom(w)\to\dom(v')$. 
\end{prop}
Finally, let us consider $N(G)$.
\begin{prop} The Neeman map
$N(G)$ (and so also $N'(G)$) is distinguished if and only if $G$ is distinguished.
\end{prop}  
\begin{proof}
  Observe that we have two 3-triangles:
  \[\xymatrix@=1.8pc{
          a'\oplus a\ar@{=}[r] &a'\oplus a \ar[r]_(0.4){\smat{1&0\\0&u\\0&0}} \ar[d]
          &a'\oplus b\oplus b'\ar[d]^(0.4){\smat{0&0&1\\0&v&0}}\\
     & 0\ar@{~>}[ul] \ar[r] & b'\oplus c \ar[r]^{\smat{0&0\\0&w}} \ar@{=}[d] & (a'\oplus
          a)[1]\ar@{=}[d]\\
     && b'\oplus c\ar@{~>}[ul] \ar[r]^{\smat{0&0\\0&w}} & (a'\oplus a)[1] 
}\]
\[\xymatrix@=1.8pc{a'\oplus a\ar[r]^{(1 \;f)} & a'\ar[r]^{u'}\ar[d]^0
          &b'\ar[d]_(0.4){\smat{v'\\0}}\\
    & a[1]\ar@{~>}[ul]^{\smat{f\\-1}} \ar[r]^{\smat{0\\1}} &c'\oplus a[1]\ar[d]^{(1\;0)}
            \ar[r]^{\smat{w'&f\\0&-1}}
            & (a'\oplus a)[1] \ar[d]^{(1\;f)}\\
      && c'\ar@{~>}[ul]^0 \ar[r]^{w'} &a'[1]      
 }\]
 
The first is degenerate and so distinguished. We show the second is isomorphic to a completion of the base
to a distinguished 3-triangle and the
maps $a_{1,2}\to a_{1,3}$ and $a_{1,3}\to a_{2,3}$ are, say, $\mat{x\\y}$ and $(m\; n)$
respectively. Then $mv'=v'$ and $w'm=w'$ so that $m$ is an isomorphism. We also have
$w'n=0$. Also $w'x+fy=f$ and $y=1$. So $w'x=0$. Then $mx+n=0$. So if we apply
$\mat{m&n\\0&1}$ at $c'\oplus a[1]$ and equality everywhere else gives an isomorphism of this
distinguished 3-triangle to the one above.

Now extend the map of bases
$\xymatrix@R=2pc@C=2.5pc@1{
  a'\oplus a\ar@{=}[r] \ar@{=}[d] & a'\oplus a \ar[d]_{(1\;f)}
  \ar[r]_(0.4){\smat{1&0\\0&u\\0&0}}  & a'\oplus b\oplus b'\ar[d]^{(u'\;g\;1)}\\
 a'\oplus a\ar[r]^{(1\;f)} & a'\ar[r]^{u'} & b' 
}$
to a distinghuished map of the 3-triangles. This gives maps
$\mat{\alpha&\beta\\ \gamma&\kappa}:b'\oplus c\to c'\oplus a[1]$ and $(\theta\;\eta):b'\oplus
c\to c'$. Then considering the restriction to the 2-triangles we have $w'\eta=fw$, $\theta=v'$,
$\alpha=v'$ and $\beta=\eta$. Then we see that $G$ is distinguished with $h=\eta=\beta$ and
$N(G)$ is distinguished with the same $h$.  
\end{proof}

In particular, since $N(G)$ can be completed to a 3-triangle in which one of the faces is the
mapping cone of $G$, we have: 
\begin{cor}
Any distinguished map of 2-triangles is good (in the sense of Neeman).
\end{cor}
Finally, we can complete the proof of Theorem \ref{t:sum} since varying $G_1$ and $G_2$ their
lightning strike equivalence classes precisely varies $G_1\oplus G_2$ in its lightning strike
equivalnce class (as a direct sum map). Then Proposition {p:sum} shows that these classes coincide.

\section{The strong 3x3 lemma}
No paper on the fundamentals of triangulated categories would be complete without a discussion
of the 3x3 lemma. It provides a useful test of the strength of the theory and in this case we
can use distinguished maps to improve the strength of the statement. Recall that the classical
theorem says that if we are given a commuting square $gu=u'f$ in a triangulated category it can be
completed to a diagram with distinguished triangles:
\begin{equation}\label{eq:3x3}\xymatrix@=1.8pc{a\ar[r]^{u}\ar[d]^{f} & b \ar[r]^{v}\ar[d]^{g} &c\ar[r]^{w}\ar[d]^{h} &a[1]\ar[d]^{f[1]}\\
  a'\ar[r]^{u'}\ar[d]^{f'} & b'\ar[r]^{v'}\ar[d]^{g'} &c'\ar[r]^{w'}\ar[d]^{h'} &a'[1]\ar[d]^{f'[1]}\\
  a''\ar[r]^{u''}\ar[d]^{f''} & b''\ar[r]^{v''}\ar[d]^{g''}
  &c''\ar[r]^{w''}\ar[d]^{h''}\ar@{}[dr]|{-} &a''[1]\ar[d]^{f''[1]}\\
  a[1]\ar[r]^{u[1]} &b[1] \ar[r]^{v[1]} &c[1] \ar[r]^{w[1]} &a[2]}
\end{equation}
where the square marked ``$-$'' anti-commutes. The classical proof goes back to 
\cite[Prop 1.1.11]{BBD} and makes a key use of the octahedral axiom. The proof shows that we can, in
fact, start with the first two rows $(a,b,c,u,v,w)$ and $(a',b',c',u',v',w')$  and first two columns $(a,a',a'',f,g,h)$
and $(b,b',b'',f',g',h')$ as distinguished triangles and then complete this to a 3x3 diagram. Neeman
further refined this in \cite[Thm 2.3]{Neeman1991} to start with the data distinguished
triangles $(a,a',a'',f,g,h)$, $(b,b',b'',f',g',h')$ and a good morphism $(f,g,h)$ of
distinguished triangles, then this can be completed to a 3x3 diagram. It is not expected that
the resulting morphism $(u,u',u'')$ of triangles would be good and so not likely that we could
specify it.  We now strengthen this to
show that if we use distinguished morphisms then we can fix both $(f,g,h)$ and $(u,u',u'')$ and
complete this to a 3x3 diagram in a fully triangulated category:

\begin{thm}\label{t:3x3}
The following diagram where $(f,g,h)$ and $(u,u',u'')$ are distinguished morphisms of
distinguished 2-triangles
\[\xymatrix@=1.8pc{a\ar[r]^{u}\ar[d]^{f} & b \ar[r]^{v}\ar[d]^{g} &c\ar[r]^{w}\ar[d]^{h} &a[1]\ar[d]^{f[1]}\\
  a'\ar[r]^{u'}\ar[d]^{f'} & b'\ar[r]^{v'}\ar[d]^{g'} &c'\ar[r]^{w'} &a'[1]\\
  a''\ar[r]^{u''}\ar[d]^{f''} & b''\ar[d]^{g''}\\
  a[1]\ar[r]^{u[1]} &b[1] }
\]
can be completed to the 3x3 diagram (\ref{eq:3x3}).
\end{thm}
\begin{proof}
Following the proof of \cite[Prop1.1.11]{BBD}, we pick a cone $d$ on $gu=u'f$ and then this fits
into two distinguished 3-triangles:
\[\xymatrix@=2pc{a\ar[r]^u &b\ar[r]^{g}\ar[d]^{v} &b'\ar[d]^{\alpha} \\
 & c\ar@{~>}[ul]^{w}\ar[r]^{\beta} & d\ar[r]^{\gamma}\ar[d]^{\delta} &a[1]\ar[d]^{u[1]}\\
 && b''\ar@{~>}[ul]^{vg''}\ar[r]^{g''} &b[1]
}\qquad
\xymatrix@=2pc{a\ar[r]^f &a'\ar[r]^{u'}\ar[d]^{f'} &b'\ar[d]^{\alpha} \\
 & a''\ar@{~>}[ul]^{f''}\ar[r]^{\theta} & d\ar[r]^{\gamma}\ar[d]^{\eta} &a[1]\ar[d]^{u[1]}\\
 && c'\ar@{~>}[ul]^{f'w'}\ar[r]^{w'} &a'[1]
}\]
Such that $\delta\alpha=g'$, $\gamma\beta=w$, $\eta\alpha=v'$ and $\gamma\theta=f''$. Recall
that in the usual proof $h$ is then defined to be $\eta\beta$. We aim to show that $\eta\beta$ differs from $h$ by a lightning strike.
To see this consider the $3$-triangle with face given by the Neeman map $N(f,g,h)$:
 \[\xymatrix@C=3.6pc@R=2.8pc{
    a'\oplus a\ar[r]^{\smat{1&0\\0&-u\\0&0}} &a'\oplus b\oplus
    b'\ar[r]^(0.6){(u',g,1)}\ar[d]^(0.4){\smat{0&0&1\\0&v&0}} &b'\ar[d]^{\smat{v'\\0}}\\
    &b'\oplus c\ar@{~>}[ul]^{\smat{0&0\\0&-w}} \ar[r]^{\smat{v'&h\\0&-w}} &c'\oplus a[1]\ar[r]^{\smat{w'&f\\0&1}}
    \ar[d]_(0.4){\smat{w'&f[1]\\ 0&-u}} &(a'\oplus a)[1]\ar[d]\\
    &&(a'\oplus b)[1]\ar@{~>}[ul]^{\smat{u'[1]&g[1]\\0&-v[1]}}\ar[r]^(0.45){\smat{1&0\\0&-1\\u'[1]&-g[1]}}
    &(a'\oplus b\oplus b')[1]
    }\]
Consider the commuting diagram of bases:
\[\xymatrix@=1.8pc{
  a\ar[r]^{u}\ar[d]^{\smat{0\\-1}} &b\ar[r]^g\ar[d]^{\smat{0\\1\\0}} &b'\ar@{=}[d]\\
a'\oplus a\ar[r] &a'\oplus b\oplus b'\ar[r]&b'
}\]
We then complete this to a distinguished map of 3-triangles. This gives three new maps:
$\smat{x\\y}:c\to b'\oplus c$, $\smat{m\\n}:d\to c'\oplus a[1]$ and $\smat{p\\q}:b''\to
(a'\oplus b')[1]$. Considering the maps of the 2-faces we have, from $(a,b,c,u,v,w)$, $w=wy$, $yv=v$
and $xv=0$. It follows that $y$ is an isomorphism. From $(a,b',d,gu,\alpha,\gamma)$ we have
$n=-\gamma$, $w'm+fn=0$, $f\gamma=w'm$ and $m\alpha=v'$. From $(b,b',b'',g,g',g'')$ we have
$p=0$ and $q=-g''$. It now follows that we have a map of distinguished 2-triangles:
\[\xymatrix@=1.8pc{a\ar[r]^{gu}\ar[d]^f &b'\ar@{=}[d]\ar[r]^{\alpha} & d\ar[d]^m\ar[r]^\gamma &a[1]\ar[d]^{f[1]}\\
a'\ar[r]^{u'} &b'\ar[r]^{v'} &c'\ar[r]^{w'} &a'[1]}\]
But this is a summand of the map $\left(\smat{0\\-1},1,\smat{m\\-\gamma}\right)$ which is a map of faces of our
distinguished map of $3$-triangles and so is distinguished.
This follows by composing with the projection
$\left((1,-f),1,(1,0)\right)$. It now follows from Theorem \ref{t:sum} that this map is
distinguished. But $(f,1,\eta)$ is also distinguished and so there is a $\tau:a[1]\to b'$ such that $m=\eta+v'\tau\gamma$. Hence, 
\begin{align*}
\eta\beta&=h-v'\tau\gamma\beta\\
&=h-v'\tau w
\end{align*}
By symmetry, $u''$ differs from $\delta\theta$ by a lightning strike. Applying these lightning strikes to $\eta $
and $\delta$ in our choice of the original two octahedra, we deduce that we can complete the 3x3 diagram with $h$ and $u''$ as required.
\end{proof}


 \bibliographystyle{alpha}
 \bibliography{p87}

\end{document}